\documentclass[11pt]{amsart}
\usepackage{amsmath}
\usepackage{amssymb}
\usepackage{amscd}

\theoremstyle{plain}
\newtheorem{thm}{Theorem}[section]
\newtheorem{lem}[thm]{Lemma}
\newtheorem{cor}[thm]{Corollary}
\newtheorem{prop}[thm]{Proposition}

\newtheorem*{ther}{Theorem}

\theoremstyle{definition}

\def \CPb {\overline{\mathbf{CP}}^{\,2}}
\def \CP {{\mathbf{CP}}^{\,2}}
\def \R {\mathbf{R}}
\def \Z {\mathbf{Z}}
\def \Sig{\Sigma}

\def \SS {S^2\times S^2}
\def \vt {\vartheta}
\def \vp {\varphi}
\def \la {\langle}
\def \ra {\rangle}
\def \a {\alpha}
\def \b {\beta}
\def \g {\gamma}
\def \d {\delta}
\def \k {\kappa}

\def \G {\Gamma}

\def \o {\omega}

\def \z {\zeta}

\def \bd {\partial}
\def \x {\times}
\def \- {\setminus}
\def \C {\subset}
\def \ve {\varepsilon}

\def \ssw {\text{SW}}
\def \sw {\mathcal{SW}}

\def \DD {\Delta}
\def \Db {\bar{\Delta}}
\def \Ds {\DD^{\text{\it{sym}}}}
\def \Dbs {\bar{\Delta}^{\text{\it{sym}}}}

\def \h {\tilde{h}}
\def \X {\tilde{X}}

\def \UL {\!|\!^{^{\overline{\;\;}}}}
\def \LL {\!|\!_{\underline{\;\;}}}
\def \UR {\;{}^{^{\overline{\;\;}}}\!|}
\def \LR {\;{}_{\underline{\;\;}}\!|}
\def \sp {\hspace{.125in}}

\begin{document}

\baselineskip.525cm

\title{Symplectic surfaces in a fixed homology class}
\author[Ronald Fintushel]{Ronald Fintushel}
\address{Department of Mathematics, Michigan State University \newline
\hspace*{.375in}East Lansing, Michigan 48824}
\email{\rm{ronfint@math.msu.edu}}
\thanks{The first author was partially supported NSF Grant DMS9704927 and the
second author by NSF Grant DMS9626330}
\author[Ronald J. Stern]{Ronald J. Stern}
\address{Department of Mathematics, University of California \newline
\hspace*{.375in}Irvine,  California 92697}
\email{\rm{rstern@math.uci.edu}}

\maketitle

\section{Introduction\label{Intro}}

The purpose of this paper is to investigate the following problem:
\begin{quotation}
For a fixed 2-dimensional homology class $\a$ in a simply connected symplectic
$4$-manifold, up to smooth isotopy, how many connected smoothly embedded
symplectic submanifolds represent $\a$?
\end{quotation}
It has been conjectured in some quarters that such a homology class $\a$ should
be represented by at most finitely many connected embedded symplectic
submanifolds; some have conjectured that such a representative must be unique.

As motivation for this
conjecture, suppose one fixes a homology class
$\a
\in H_2(X;\Z)$ where $X$
is a Kahler surface and asks, up to smooth isotopy, how many nonsingular
complex
curves represent this class. If
$a\in H^2(X;\Z)$ is the Poincar\'e dual of $\a$, then each complex curve
representing
$\a$ is the zero set of a section of a holomorphic line bundle with
$c_1=a$.  Thus
we must ask about the preimage of $a$ under the map
\[ c_1:H^1(X;\mathcal{O}^*_X)\to H^2(X;\Z).\]
Equivalently, we study the kernel of
$c_1$. This is an analytic variety, and hence has finitely many connected
components.
Since the points corresponding to singular curves form a subvariety of complex
codimension at least one, up to smooth isotopy, there are at most
finitely many nonsingular
complex curves representing $\a$. In particular, if $H^1(X;\Z)=0$ then the
Picard
torus, $\ker(c_1)=0$; so there is a unique representative of the class in
question.

In contrast we shall prove the following theorem in \S{5}.

\begin{ther} Let $X$ be a simply connected symplectic $4$-manifold which
contains a
c-embedded symplectic torus $T$. Then in each homology class $2m\,[T]$,
$m\ge2$,
there is an infinite family of smoothly embedded symplectic tori, no two of
which
are smoothly isotopic.
\end{ther}

\noindent To say that a torus $T$ is {\it{c-embedded}} means that $T$ is a
smoothly
embedded homologically essential torus of self-intersection zero which has
a a pair
of simple curves which generate its first homology and which bound
vanishing cycles
(disks of self-intersection $-1$) in $X$. (See \cite{KL4M}.) The simplest
examples
of c-embedded tori are generic fibers of simply connected elliptic
fibrations. One
can also find c-embedded tori in many surfaces of general type (including
Horikawa
surfaces) via the process of rational blowdowns \cite{rat}.

One might then ask for what families of symplectic $4$-manifolds are there only
finitely many smooth isotopy classes of symplectic surfaces in any fixed
homology class. In light of the above theorem,
a reasonable conjecture might be that this finiteness condition
holds for
ruled surfaces or rational surfaces with $c_1^2>0$. Siebert and Tian have shown
that each symplectic surface in $\SS$ with genus $\le3$ is smoothly
isotopic to a
complex curve; so finiteness holds in that situation.

The technique of this paper, described in detail below, is to replace the
torus $T$
in its tubular neighborhood, $T=S^1\x S^1\C S^1\x S^1\x D^2$, with
$S^1\x B\C S^1\x S^1\x D^2$ where $B$ is a closed braid in $S^1\x D^2$. In
case the
braid $B$ has an even number of strands and also represents the unknot in
$S^3$, let  $L_B$ denote the 2-component link in $S^3$
obtained as the preimage of the axis of $B$ under the $2$-fold cover of $S^3$
branched over $B$. We can then
identify the double cover of the symplectic manifold $X$ branched over
$S^1\x B$
as the manifold $X_{L_B}$ of \cite{KL4M}.  This manifold has a
Seiberg-Witten invariant which was
computed in
\cite{KL4M}; it is related to the Alexander polynomial of the link $L_B$.
Note that if $B$ has $2m$ strands,
$S^1\x B$ is homologous to
$2mT$. In order to obtain infinitely many nonisotopic such tori homologous
to $2mT$,  we will
utilize a braid construction of Birman and Menasco \cite{BM5,BM6} to
construct infinitely many such braids
which are distinguished by their Alexander polynomials.

The construction of this paper contrasts with an older construction of the
authors that produced (under mild hypotheses) infinitely many non-smoothly
isotopic
embedded surfaces, all topologically ambiently isotopic to a given embedded
surface
\cite{surfaces}. This older
construction replaced an annulus $S^1\x I\x\{ 0\}\C S^1\x I\x D^2$ with
$S^1\x K\C S^1\x I\x D^2$ where $K$ is the result of tying a knot in the core
$I\x\{ 0\}$ of the cylinder $I\x D^2$. If $\Sig$ is a symplectic surface of
positive
genus and nonnegative self-intersection, and $\Sig_K$ is the result of
performing
this knotting operation, it is shown in \cite{surfaces} that $\Sig_K$ is
{\it not}
smoothly isotopic to a symplectic submanifold as long as the Alexander
polynomial of
$K$ is nontrivial.

There have also been informal conjectures asserting the finiteness of the
number
of smooth simply connected symplectic $4$-manifolds in a fixed
homeomorphism type
which admit symplectic Lefschetz fibrations with a fiber of fixed genus. In the
last section of this paper we produce counterexamples, stemming from our
examples,
to these conjectures. Other, more easily obtained, counterexamples are given
in \cite{FSLef}.

It is interesting to ask whether the nonfiniteness results of this paper are a
general phenomenon applying to surfaces of arbitrary genus or whether they are
unique to tori. The authors have general constructions which apply to surfaces
of higher genus, but they have been unable to determine whether or not the
resulting surfaces are smoothly isotopic.

Finally, the authors wish to express appreciation to Bill Menasco for
(e-mail) conversations concerning his joint work with Joan Birman, and to
Gang Tian
whose interest stimulated this work.

\section{Braids}\label{braids}

In this section we shall describe a sequence of families of closed $2m$-strand
braids
$B_{2m,k}$ in $S^3$, $m,k=1,2,\dots$, whose corresponding double branched
covers
yield fibered 2-component links in $S^3$.
We begin by describing a construction of fibered 2-component links due to D.
Goldsmith
\cite{G}. Let $B$ be a closed $2m$-strand braid in $S^3$ with axis $A$.
I.e., $B$
is a braid in an unknotted solid torus $V=S^1\x D^2$ in $S^3$, and $A$ is
the core
of the complementary unknotted solid torus. We may think of $A$ as a
fibered knot,
whose fibers are the disks $\{ t_0\}\x D^2$ of $V$. Each such disk contains
$2m$
points of the braid $B$.

\vspace{-.5in}

\centerline{\unitlength 1.5cm
\begin{picture}(9,6)
\qbezier(2.75,4)(4.5,5)(6.25,4)
\qbezier(2.75,2)(4.5,1)(6.25,2)
\qbezier(2.75,4)(1,3)(2.75,2)
\qbezier(6.25,4)(8,3)(6.25,2)
\qbezier(2.85,3.1)(3.25,2.35)(3.65,3.1)
\qbezier(2.9,3)(3.25,3.6)(3.6,3)
\qbezier(5.35,3.1)(5.75,2.35)(6.15,3.1)
\qbezier(5.4,3)(5.75,3.6)(6.1,3)
\qbezier(2.9,3)(2.4,3.735)(1.875,3)
\qbezier[20](2.9,3)(2.4,2.265)(1.875,3)
\put (2.4,3.45){\SMALL{$C_1$}}
\put(3.25,3.05){\oval(1.15,.85)}
\put(3.15,3.6){\SMALL{$C_2$}}
\qbezier(4,3.425)(3.65,3.32)(3.6,3)
\qbezier[10](4,2.575)(3.65,2.68)(3.6,3)
\put(3.85,3.6){\SMALL{$C_3$}}
\put(4.25,3){\LARGE${\bf{\dots}}$}
\qbezier(7.125,3)(6.6,3.735)(6.075,3)
\qbezier[20](7.125,3)(6.6,2.265)(6.075,3)
\put (6.3,3.45){\SMALL{$C_{2m-1}$}}
\put(5.75,3.05){\oval(1.15,.85)}
\put(5.5,3.6){\SMALL{$C_{2m-2}$}}
\qbezier(5,3.425)(5.35,3.32)(5.4,3)
\qbezier[10](5,2.575)(5.35,2.68)(5.4,3)
\put(4.5,4){\circle*{.5}}
\put(4.5,2){\circle*{.5}}
\put(6.5,1.75){$\Sig_{m-1}^{''}$}
\put(4.15,1) {Figure 1}
\end{picture}}

\vspace{-.25in}
Now suppose further that $B$ represents an unknotted circle in $S^3$, that
is to
say, $B$ can be isotoped to the unknot in $S^3$ when one allows it to pass
through
$A$. The double branched cover of $S^3$ branched over $B$ is then $S^3$
again, and
since $A$ links $B$ an even number of times, it lifts to a 2-component link
$L_B$ in
the cover $S^3$. This link is fibered, and its fibers are simply the  double
branched covers of the fibers of the unknot $A$. These are twice-punctured
surfaces
of genus $m-1$. So we have
\[ S^3\- L_B = S^1\x_\vp \Sig''_{m-1}\] where $\Sig''_{m-1}$ is the surface of
genus $m-1$ with two boundary components. The monodromy map $\vp$ can be
calculated
from the braid $B$. The braid group on $2m$ strands is generated by the
elementary
braid transpositions $\b_1,\dots,\b_{2m-1}$, where $\b_i$ denotes a right-hand
crossing of the $i$th strand over the $(i+1)$st. In the double branched
cover each
such crossing contributes a Dehn twist. (See e.g.\cite[p.172]{BZ}.)  If we
write the
braid group element corresponding to $B$ as a word in the $\{\b_i\}$, it
follows
that the monodromy will be the product of Dehn twists about the simple
closed
curves $\{ C_i\}$ as shown in Figure~1.

Our next task is to construct for each integer $n \ge 4$ a family of closed
n-braids $\{ B_{n,k}\}$,
$k$ a nonnegative integer, with the properties that each $B_{n,k}$ is
unknotted in $S^3$ and, for fixed $m$,
the Alexander polynomial of the 2-component link $L_{B_{2m,k}}$ is
distinguished by the integer $k$. It is these
braids  which will be used to construct our examples of symplectic
submanifolds.

We begin with the 4-strand braid
$B_{4,0}$, shown in Figure~2, first constructed by Birman and Menasco
\cite{BM5}. For us, the key property of this braid is that it represents
the unknot in
$S^3$.

\centerline{\unitlength .8cm
\begin{picture}(7,8)
\put (1,7.5){\line(1,0){5}}
\put (1,7.5){\line(0,-1){4}}
\put (1.5,7){\line(1,0){4}}
\put (1.5,7){\line(0,-1){2.5}}
\put (2,6.5){\line(1,0){1.25}}
\put (2,6.5){\line(0,-1){2}}
\qbezier (1.5,4.5)(1.75,4.25)(2,4)
\qbezier (2,4.5)(1.9,4.4)(1.8,4.3)
\qbezier (1.5,4)(1.6,4.1)(1.7,4.2)
\qbezier (1.5,4)(1.75,3.75)(2,3.5)
\qbezier (2,4)(1.9,3.9)(1.8,3.8)
\qbezier (1.5,3.5)(1.6,3.6)(1.7,3.7)
\qbezier (1.5,3.5)(1.25,3.25)(1,3)
\qbezier (1,3.5)(1.1,3.4)(1.2,3.3)
\qbezier (1.5,3)(1.4,3.1)(1.3,3.2)
\put (2,3.5){\line(0,-1){.5}}
\qbezier (1.5,3)(1.75,2.75)(2,2.5)
\qbezier (2,3)(1.9,2.9)(1.8,2.8)
\qbezier (1.5,2.5)(1.6,2.6)(1.7,2.7)
\put (2.5,6){\line(0,-1){4.25}}
\put (2.5,6){\line(1,0){.75}}
\qbezier (3.25,6)(3.5,6.25)(3.75,6.5)
\qbezier (3.25,6.5)(3.35,6.4)(3.45,6.3)
\qbezier (3.75,6)(3.65,6.1)(3.55,6.2)
\put (3.75,6){\line(1,0){.75}}
\put (3.75,6.5){\line(1,0){1.25}}
\put (5.5,7) {\line(0,-1){1.75}}
\put (5,6.5) {\line(0,-1){1.25}}
\qbezier (5.5,5.25)(5.25,5)(5,4.75)
\qbezier (5,5.25)(5.1,5.15)(5.2,5.05)
\qbezier (5.5,4.75)(5.4,4.85)(5.3,4.95)
\put (6,7.5) {\line(0,-1){2.75}}
\qbezier (5.5,4.75)(5.75,4.5)(6,4.25)
\qbezier (6,4.75)(5.9,4.65)(5.8,4.55)
\qbezier (5.5,4.25)(5.6,4.35)(5.7,4.45)
\put (5,4.75) {\line(0,-1){.5}}
\qbezier (5.5,4.25)(5.25,4)(5,3.75)
\qbezier (5,4.25)(5.1,4.15)(5.2,4.05)
\qbezier (5.5,3.75)(5.4,3.85)(5.3,3.95)
\qbezier (5.5,3.75)(5.25,3.5)(5,3.25)
\qbezier (5,3.75)(5.1,3.65)(5.2,3.55)
\qbezier (5.5,3.25)(5.4,3.35)(5.3,3.45)
\qbezier (5.5,3.25)(5.25,3)(5,2.75)
\qbezier (5,3.25)(5.1,3.15)(5.2,3.05)
\qbezier (5.5,2.75)(5.4,2.85)(5.3,2.95)
\put (4.5,6) {\line(0,-1){4.25}}
\put (5,2.75) {\line(0,-1){1.5}}
\put (5.5,2.75) {\line(0,-1){2}}
\put (6,4.25) {\line(0,-1){4}}
\put (2,2.5) {\line(0,-1){1.25}}
\put (1.5,2.5) {\line(0,-1){1.75}}
\put (1,3) {\line(0,-1){2.75}}
\put (1,.25) {\line(1,0){5}}
\put (1.5,.75) {\line(1,0){4}}
\qbezier (3.25,1.75)(3.5,1.5)(3.75,1.25)
\qbezier (3.25,1.25)(3.35,1.35)(3.45,1.45)
\qbezier (3.75,1.75)(3.65,1.65)(3.55,1.55)
\put (2.5,1.75){\line(1,0){.75}}
\put (3.75,1.75){\line(1,0){.75}}
\put (2,1.25){\line(1,0){1.25}}
\put (3.75,1.25){\line(1,0){1.25}}
\put (2.5,-.5){Figure 2:  $B_{4,0}$}
\put(3.5,4){\LARGE{$\bf{.}$}}
\put(3.5,3.75){$A$}
\end{picture}}

\vspace{.25in}

Using  the integer $j$ as shorthand for the braid transposition $\b_j$ and
$\bar{j}$ for $\b_j^{-1}$, the braid $B_{4,0}$ is given by the expression
\[ B_{4,0} = (\bar{2}.\bar{2}.1.\bar{2}).3.(2.2.2.\bar{1}.2).\bar{3}\]
We define braids $B_{m,0}$ inductively as follows. Assume that $B_{m,0}$ is
given
by
\[B_{m,0}=\Phi_m.(m-1).\Psi_m.(\overline{m-1})\]
where $\Phi_m$ and $\Psi_m$ are expressions involving only the braid
transpositions
$j<m-1$ and their inverses. Define
\[B_{m+1,0}=(m-2).(m-1).\Phi_m.(\overline{m-1}).m.(\overline{m-1}).
\Psi_m.(m-1).\bar{m}\]
Thus $\Phi_{m+1}=(m-2).(m-1).\Phi_m.(\overline{m-1})$ and
$\Psi_{m+1}=(\overline{m-1}).\Psi_m.(m-1)$.
A schematic is given in Figure~3.

\bigskip
\centerline{\unitlength .8cm
\begin{picture}(14,8)
\put (1,7.5){\line(1,0){5}}
\put (1,7.5){\line(0,-1){7.25}}
\put (1.5,7){\line(1,0){4}}
\put (1.5,7){\line(0,-1){2.5}}
\put (1,3){\line(1,0){1}}
\put (1,4.5){\line(1,0){1}}
\put (2,6.5){\line(1,0){1.25}}
\put (2,6.5){\line(0,-1){5.25}}
\put (2.5,6){\line(0,-1){4.25}}
\put (2.5,6){\line(1,0){.75}}
\put (3.75,6){\line(1,0){.75}}
\put (3.75,6.5){\line(1,0){1.25}}
\put (5.5,7) {\line(0,-1){2.5}}
\put (5,6.5) {\line(0,-1){5.25}}
\put (6,7.5) {\line(0,-1){7.25}}
\put (5,4.75) {\line(0,-1){.5}}
\put (4.5,6) {\line(0,-1){4.25}}
\put (5.5,3) {\line(0,-1){2.25}}
\put (1.5,3) {\line(0,-1){2.25}}
\put (5,3){\line(1,0){1}}
\put (5,4.5){\line(1,0){1}}
\put (1,.25) {\line(1,0){5}}
\put (1.5,.75) {\line(1,0){4}}
\put (2.5,1.75){\line(1,0){.75}}
\put (3.75,1.75){\line(1,0){.75}}
\put (2,1.25){\line(1,0){1.25}}
\put (3.75,1.25){\line(1,0){1.25}}
\qbezier (3.25,6)(3.5,6.25)(3.75,6.5)
\qbezier (3.25,6.5)(3.35,6.4)(3.45,6.3)
\qbezier (3.75,6)(3.65,6.1)(3.55,6.2)
\qbezier (3.25,1.75)(3.5,1.5)(3.75,1.25)
\qbezier (3.25,1.25)(3.35,1.35)(3.45,1.45)
\qbezier (3.75,1.75)(3.65,1.65)(3.55,1.55)
\put (1.2,3.6){\!\Large{$\Phi_m$}}
\put (5.2,3.6){\!\Large{$\Psi_m$}}
\put (3.15,-.25){$B_{m,0}$}
\put(3.4,4){\LARGE{$\bf{.}$}}
\put(3.4,3.75){$A$}
\put (8,7.5){\line(1,0){5}}
\put (8,7.5){\line(0,-1){7.25}}
\put (8.5,7.25){\line(1,0){4}}
\put (8.5,7.25){\line(0,-1){2.75}}
\put (8,3){\line(1,0){1}}
\put (8,4.5){\line(1,0){1}}
\put (9,6.5){\line(1,0){1.25}}
\put (9,6.5){\line(0,-1){5.25}}
\put (9.5,6){\line(0,-1){4.25}}
\put (9.5,6){\line(1,0){.75}}
\put (10.75,6){\line(1,0){.75}}
\put (10.75,6.5){\line(1,0){.9}}
\put (12,6.5){\line(-1,0){.15}}
\put (12.5,7.25) {\line(0,-1){2.75}}
\put (12,6.5) {\line(0,-1){5.25}}
\put (13,7.5) {\line(0,-1){7.25}}
\put (12,4.75) {\line(0,-1){.5}}
\put (11.5,6) {\line(0,-1){4.25}}
\put (12.5,3) {\line(0,-1){2.25}}
\put (8.5,3) {\line(0,-1){2.25}}
\put (12,3){\line(1,0){1}}
\put (12,4.5){\line(1,0){1}}
\put (8,.25) {\line(1,0){5}}
\put (8.5,.75) {\line(1,0){4}}
\put (9.5,1.75){\line(1,0){.75}}
\put (10.75,1.75){\line(1,0){.75}}
\put (9,1.25){\line(1,0){1.25}}
\put (10.75,1.25){\line(1,0){.9}}
\put (12,1.25){\line(-1,0){.15}}
\qbezier (10.25,6)(10.5,6.25)(10.75,6.5)
\qbezier (10.25,6.5)(10.35,6.4)(10.45,6.3)
\qbezier (10.75,6)(10.65,6.1)(10.55,6.2)
\qbezier (10.25,1.75)(10.5,1.5)(10.75,1.25)
\qbezier (10.25,1.25)(10.35,1.35)(10.45,1.45)
\qbezier (10.75,1.75)(10.65,1.65)(10.55,1.55)
\put (8.2,3.6){\!\Large{$\Phi_m$}}
\put (12.2,3.6){\!\Large{$\Psi_m$}}
\put (8.75,4.5) {\line(0,1){2.25}}
\put (8.75,6.75) {\line(1,0){3}}
\put (11.75,6.75) {\line(0,-1){5.75}}
\put (11.75,1) {\line(-1,0){2.5}}
\put (9.25,1) {\line(0,1){.15}}
\put (9.25,1.35) {\line(0,1){5.1}}
\put (9.25,6.6) {\line(0,1){.15}}
\put (9.25,6.85) {\line(0,1){.15}}
\put (9.25,7) {\line(1,0){3}}
\put (12.25,7) {\line(0,-1){2.5}}
\put (10.15,-.25){$B_{m+1,0}$}
\put(10.4,4){\LARGE{$\bf{.}$}}
\put(10.4,3.75){$A$}
\put (6.5,-.75){Figure 3}
\end{picture}}

\vspace{.25in}

\begin{lem}\label{unknotted} When the braids $B_{m,0}$ ($m\ge4$) are
considered as
knots in $S^3$, they are unknotted.\end{lem}
\begin{proof} Figure~4 shows how $B_{m+1,0}$ is isotopic to $B_{m,0}$ in
$S^3$, and
this completes the proof since $B_{4,0}$ is unknotted.
\end{proof}

\bigskip
\centerline{\unitlength .8cm
\begin{picture}(14,8)
\put (1,7.5){\line(1,0){5}}
\put (1,7.5){\line(0,-1){7.25}}
\put (1.5,7.25){\line(1,0){4}}
\put (1.5,7.25){\line(0,-1){2.75}}
\put (1,3){\line(1,0){1}}
\put (1,4.5){\line(1,0){1}}
\put (2,6.5){\line(1,0){1}}
\put (2,6.5){\line(0,-1){5.25}}
\put (5,6.5){\line(-1,0){1}}
\put (5.5,7.25) {\line(0,-1){2.75}}
\put (5,6.5) {\line(0,-1){5.25}}
\put (6,7.5) {\line(0,-1){7.25}}
\put (5,4.75) {\line(0,-1){.5}}
\put (5.5,3) {\line(0,-1){2.25}}
\put (1.5,3) {\line(0,-1){2.25}}
\put (5,3){\line(1,0){1}}
\put (5,4.5){\line(1,0){1}}
\put (1,.25) {\line(1,0){5}}
\put (1.5,.75) {\line(1,0){4}}
\put (2,1.25){\line(1,0){1}}
\put (3,1.25){\line(0,1){5.25}}
\put (5,1.25){\line(-1,0){1}}
\put (4,1.25){\line(0,1){5.25}}
\put (1.2,3.6){\!\Large{$\Phi_m$}}
\put (5.2,3.6){\!\Large{$\Psi_m$}}
\put (1.75,4.5) {\line(0,1){2.25}}
\put (1.75,6.75) {\line(1,0){1.75}}%
\put (2.25,1) {\line(1,0){1.25}}%
\put (3.5,1) {\line(0,1){5.725}}%
\put (2.25,1) {\line(0,1){.15}}
\put (2.25,1.35) {\line(0,1){5.1}}
\put (2.25,6.6) {\line(0,1){.15}}
\put (2.25,6.85) {\line(0,1){.15}}
\put (2.25,7) {\line(1,0){3}}
\put (5.25,7) {\line(0,-1){2.5}}
\put (3,-.25){`$B_{m+1,0}$'}
\put (8,7.5){\line(1,0){5}}
\put (8,7.5){\line(0,-1){7.25}}
\put (8.5,7.25){\line(1,0){4}}
\put (8.5,7.25){\line(0,-1){2.75}}
\put (8,3){\line(1,0){1}}
\put (8,4.5){\line(1,0){1}}
\put (9,6.5){\line(1,0){1}}
\put (9,6.5){\line(0,-1){5.25}}
\put (12,6.5){\line(-1,0){1}}
\put (12.5,7.25) {\line(0,-1){2.75}}
\put (12,6.5) {\line(0,-1){5.25}}
\put (13,7.5) {\line(0,-1){7.25}}
\put (12,4.75) {\line(0,-1){.5}}
\put (12.5,3) {\line(0,-1){2.25}}
\put (12.5,3) {\line(0,-1){2.25}}
\put (12,3){\line(1,0){1}}
\put (12,4.5){\line(1,0){1}}
\put (8,.25) {\line(1,0){5}}
\put (8.5,.75) {\line(1,0){4}}
\put (9,1.25){\line(1,0){1}}
\put (10,1.25){\line(0,1){5.25}}
\put (12,1.25){\line(-1,0){1}}
\put (11,1.25){\line(0,1){5.25}}
\put (8.2,3.6){\!\Large{$\Phi_m$}}
\put (12.2,3.6){\!\Large{$\Psi_m$}}
\put (8.75,4.5) {\line(0,1){2.5}}
\put (8.75,7) {\line(1,0){3.5}}
\put (12.25,7) {\line(0,-1){2.5}}
\put (8.5,3) {\line(0,-1){2.25}}
\put (10,-.25){`$B_{m,0}$'}
\put (6.5,-.75){Figure 4}
\end{picture}}

\vspace{.5in}

In \cite{BM6}, Birman and Menasco introduced an operation on m-strand
braids of the
form $B= \Phi.(m-1).\Psi.(\overline{m-1})$ where where $\Phi$ and $\Psi$ are
expressions in braid transpositions $j<m-1$ and their inverses. This
operation is
pictured in Figure~5.

\bigskip

\centerline{\unitlength .8cm
\begin{picture}(14,8)
\put (1,7.5){\line(1,0){5}}
\put (1,7.5){\line(0,-1){7.25}}
\put (1.5,7){\line(1,0){4}}
\put (1.5,7){\line(0,-1){2.5}}
\put (1,3){\line(1,0){1}}
\put (1,4.5){\line(1,0){1}}
\put (2,6.5){\line(1,0){1.25}}
\put (2,6.5){\line(0,-1){5.25}}
\put (2.5,6){\line(0,-1){4.25}}
\put (2.5,6){\line(1,0){.75}}
\put (3.75,6){\line(1,0){.75}}
\put (3.75,6.5){\line(1,0){1.25}}
\put (5.5,7) {\line(0,-1){2.5}}
\put (5,6.5) {\line(0,-1){5.25}}
\put (6,7.5) {\line(0,-1){7.25}}
\put (5,4.75) {\line(0,-1){.5}}
\put (4.5,6) {\line(0,-1){4.25}}
\put (5.5,3) {\line(0,-1){2.25}}
\put (1.5,3) {\line(0,-1){2.25}}
\put (5,3){\line(1,0){1}}
\put (5,4.5){\line(1,0){1}}
\put (1,.25) {\line(1,0){5}}
\put (1.5,.75) {\line(1,0){4}}
\put (2.5,1.75){\line(1,0){.75}}
\put (3.75,1.75){\line(1,0){.75}}
\put (2,1.25){\line(1,0){1.25}}
\put (3.75,1.25){\line(1,0){1.25}}
\qbezier (3.25,6)(3.5,6.25)(3.75,6.5)
\qbezier (3.25,6.5)(3.35,6.4)(3.45,6.3)
\qbezier (3.75,6)(3.65,6.1)(3.55,6.2)
\qbezier (3.25,1.75)(3.5,1.5)(3.75,1.25)
\qbezier (3.25,1.25)(3.35,1.35)(3.45,1.45)
\qbezier (3.75,1.75)(3.65,1.65)(3.55,1.55)
\put (1.3,3.6){\Large{$\Phi$}}
\put (5.3,3.6){\Large{$\Psi$}}
\put(3.4,4){\LARGE{$\bf{.}$}}
\put(3.4,3.75){$A$}
\put (6.75,3.75){$\longrightarrow$}
\put (8,7.5){\line(1,0){5}}
\put (8,7.5){\line(0,-1){1.35}}
\put (8,5.85){\line(0,-1){4}}
\put (8,.25){\line(0,1){1.35}}
\put (8.5,7){\line(1,0){4}}
\put (8.5,7){\line(0,-1){.6}}
\put (8.5,6.15){\line(0,-1){1.65}}
\put (8,3){\line(1,0){1}}
\put (8,4.5){\line(1,0){1}}
\put (9,6.5){\line(1,0){1.25}}
\qbezier(9,6.5)(8.25,6.15)(7.9,5.95)
\qbezier(7.9,5.95)(7.5,5.75)(7.9,5.55)
\qbezier(8.1,5,45)(8.25,5.375)(8.4,5.3)
\qbezier(8.6,5.2)(8.75,5.125)(9,5)
\put (9,5){\line(0,-1){2.25}}
\qbezier(7.9,2.2)(7.5,2)(7.9,1.8)
\qbezier(9,1.25)(8.45,1.5025)(7.9,1.8)
\qbezier(9,2.75)(8.75,2.625)(8.6,2.55)
\qbezier(8.4,2.45)(8.25,2.375)(8.1,2.3)
\put (9.5,6){\line(0,-1){4.25}}
\put (9.5,6){\line(1,0){.75}}
\put (10.75,6){\line(1,0){.75}}
\put (10.75,6.5){\line(1,0){1.25}}
\put (12.5,7) {\line(0,-1){2.5}}
\put (12,6.5) {\line(0,-1){5.25}}
\put (13,7.5) {\line(0,-1){7.25}}
\put (12,4.75) {\line(0,-1){.5}}
\put (11.5,6) {\line(0,-1){4.25}}
\put (12.5,3) {\line(0,-1){2.25}}
\put (8.5,3) {\line(0,-1){1.35}}
\put (8.5,.75) {\line(0,1){.65}}
\put (12,3){\line(1,0){1}}
\put (12,4.5){\line(1,0){1}}
\put (8,.25) {\line(1,0){5}}
\put (8.5,.75) {\line(1,0){4}}
\put (9.5,1.75){\line(1,0){.75}}
\put (10.75,1.75){\line(1,0){.75}}
\put (9,1.25){\line(1,0){1.25}}
\put (10.75,1.25){\line(1,0){1.25}}
\qbezier (10.25,6)(10.5,6.25)(10.75,6.5)
\qbezier (10.25,6.5)(10.35,6.4)(10.45,6.3)
\qbezier (10.75,6)(10.65,6.1)(10.55,6.2)
\qbezier (10.25,1.75)(10.5,1.5)(10.75,1.25)
\qbezier (10.25,1.25)(10.35,1.35)(10.45,1.45)
\qbezier (10.75,1.75)(10.65,1.65)(10.55,1.55)
\put (8.3,3.6){\Large{$\Phi$}}
\put (12.3,3.6){\Large{$\Psi$}}
\put(10.4,4){\LARGE{$\bf{.}$}}
\put(10.4,3.75){$A$}
\put (6.5,-.75){Figure 5}
\end{picture}}

\vspace{.5in}

The Birman-Menasco operation preserves the link-type of the braid (as a link in
$S^3$). Formally, the Birman-Menasco operation is:
\[ \Phi.(m-1).\Psi.(\overline{m-1}) \;\;\longrightarrow \;\;
\G_{m-2}^{-1}.\Phi.\G_{m-2}.(m-1).\Psi.(\overline{m-1}) \]
where
\[\G_r=r.(r-1).\cdots.2.1.1.2.\cdots.(r-1).r \]

Our family of m-strand braids is $\{ B_{m,k}\}$ where $B_{m,k}$ is the
result of
applying the Birman-Menasco operation $k$ times to the braid $B_{m,0}$.
Hence
\[ B_{m,k} = \G_{m-2}^{-k}.\Phi_m.\G_{m-2}^k.(m-1).\Psi_m.(\overline{m-1})\]
or $B_{m,k} = \Phi_{m,k}.(m-1).\Psi_m.(\overline{m-1})$ where
$\Phi_{m,k}=\G_{m-2}^{-k}.\Phi_m.\G_{m-2}^k$. It follows from
Lemma~\ref{unknotted},
and the fact (easily seen in Figure~5) that the Birman-Menasco operation
preserves
the link type of the braid, that the braids $B_{m,k}$
all represent unknots in
$S^3$.

\section{The double covering links}

In this section we shall study the 2-component links $L_{2m,k}$ which
result from
taking the preimage \,$\pi^{-1}(A)$\, of the axis in the double cover of $S^3$
branched over $B_{2m,k}$. Recall from \S~\!\!2 that $L_{2m,k}$ is a fibered
link and
its fiber is the twice-punctured surface $\Sig_{m-1}''$ of genus $m-1$. We are
interested in the monodromy of this fibration. As was discussed in
\S~\!\!2, this monodromy is a product of Dehn twists given by the braid
transpositions which describe $B_{2m,k}$ as an element of the braid group on
$2m$-strands. Each transposition $\b_j$ corresponds to the Dehn twist about the
curve $C_j$ of Figure~1. We orient these curves so that their intersection
numbers
are
  \begin{eqnarray*}C_i\cdot C_j& = &0,\; \;j\ne i\pm 1\\ C_{i-1}\cdot
C_i&=&1 \\
C_i\cdot C_{i+1}&=&-1 \end{eqnarray*}
In homology, the Dehn twist corresponding to $\b_k$ is given by
$a\to a+(a\cdot C_k)\,C_k$. Thus the matrix representing this Dehn twist on
$H_1(\Sig_{m-1}'';\Z)$ is $D_{2m,k}=I_{2m-1}+J_{2m-1,k}$ where $I_{2m-1}$
is the
identity matrix of rank $2m-1$ and $J_{2m-1,k}$ is the $(2m-1)\x(2m-1)$ matrix
whose entries are all $0$ except for $(J_{2m-1,k})_{{}_{k,k-1}}=1$ and
$(J_{2m-1,k})_{{}_{k,k+1}}=-1$.

Denote the matrix representing the homology monodromy of $L_{2m,k}$
by $\Omega_{2m,k}$. For example, $B_{4,0} =
(\bar{2}.\bar{2}.1.\bar{2}).3.(2.2.2.\bar{1}.2).\bar{3}$; so
\[ \Omega_{4,0}= D_{4,2}^{-1}\cdot D_{4,2}^{-1}\cdot D_{4,1}\cdot
D_{4,2}^{-1}\cdot
     D_{4,3}\cdot  D_{4,2}\cdot D_{4,2}\cdot D_{4,2}\cdot D_{4,1}^{-1}\cdot
     D_{4,2}\cdot D_{4,3}^{-1} \]
This matrix is
\[ \Omega_{4,0}=\begin{pmatrix} -10 & -17 & \sp 11\\ \sp 46 & \sp 73 & -46\\
\sp 7 & \sp 10 &\hspace{.05in} -6
\end{pmatrix} \]

In order to save notation we shall denote by $\Phi_n$, $\Phi_{n,k}$, and
$\Psi_n$
the rank $n-1$ square matrices corresponding to the product of Dehn twists
resulting
from the braid group elements with the same name. (For these purely
combinatorial
expressions, there is no need to assume that $n$ is even.) Similarly, we let
$\G_{n,r}$ be the rank
$n-1$ square matrix corresponding to $\G_r$. An easy inductive argument gives:

\begin{lem}\label{Gamma} For any integer $k$, the matrix power
$\G_{n,n-2}^{\,k}$ is given by:
\medskip
\[\begin{pmatrix} {}&{}&{}&2k&0\\
{}&{}&{}&0&0\\
{}&{}&{}&2k&0\\
{}&{}&{}&0&0\\
{}&\parbox[r]{.75in}{\hspace*{.25in}\LARGE{$I_{n-3}$}}&{}&\vdots&\vdots\\
{}&{}&{}&2k&0\\
{}&{}&{}&0&0\\
{}&{}&{}&2k&0\\
0&\cdots&0&1&0\\
0&\cdots&0&0&1
\end{pmatrix} \;\text{$k$ even},
\hspace{.25in}
\begin{pmatrix} {}&{}&{}&-2k&2\\
{}&{}&{}&\sp 0&0\\
{}&{}&{}&-2k&2\\
{}&{}&{}&\sp 0&0\\
{}&\parbox[r]{.75in}{\hspace*{.15in}\LARGE{$-I_{n-3}$}}&{}&\sp\vdots&\vdots\\
{}&{}&{}&-2k&2\\
{}&{}&{}&\sp 0&0\\
{}&{}&{}&-2k&2\\
0&\cdots&0&-1&0\\
0&\cdots&0&\sp 0&1
\end{pmatrix}
  \;\text{$k$ odd.} \]
\end{lem}
\bigskip
\bigskip
Recall that the matrices $\Phi_{n,k}$ are
recursively defined by the formulas
\begin{eqnarray*}\Phi_{n,k}&=&\G_{n,n-2}^{-k}.\Phi_{n}.\G_{n,n-2}^k\\
\Phi_{n+1}&=&(n-2).(n-1).\Phi_n.(\overline{n-1})\end{eqnarray*}

\noindent Thus

\begin{eqnarray*}\Phi_{2m,k}&=&\G_{2m,2m-2}^{-k}.\Phi_{2m}.\G_{2m,2m-2}^k\\
\Phi_{2m+2}&=&(2m-1).(2m).(2m-2).(2m-1).\Phi_{2m}.(\overline{2m-1}).(\
overline{2
m})\end{eqnarray*}

Since we have
$\Phi_4 \; = \; \begin{pmatrix}\sp 2&-1&-1\\-5&\sp 3&\sp 5\\
\sp 0&\sp 0&\sp 1\end{pmatrix}$, we obtain the following closed formulas for
$\Phi_{2m,k}$ (and consequently for $\Phi_{2m}=\Phi_{2m,0}$).

\begin{lem}\label{Phi} For $m\ge 3$, the matrices $\Phi_{2m,k}$ are
given by
\[ \left(\begin{matrix}
10k+2&6k&0&0&0&0&0&\cdots&0\\
2&1&-1&0&0&0&0&\cdots&0\\
10k+2&6k+1&0&-1&0&0&0&\cdots&0\\
2&1&0&0&-1&0&0&\cdots&0\\
10k+2&6k+1&0&0&0&-1&0&\cdots&0\\
2&1&0&0&0&0&-1&\cdots&0\\
\vdots&\vdots&\vdots&\vdots&\vdots&\vdots&\vdots&\cdots&\vdots\\
10k+2&6k+1&0&0&0&0&0&\cdots&-1\\
2&1&0&0&0&0&0&\cdots&0\\
10k+2&6k+1&0&0&0&0&0&\cdots&0\\
-5&-3&0&0&0&0&0&\cdots&0\\
0&0&0&0&0&0&0&\cdots&0 \end{matrix}
\begin{matrix}
0&20k^2-8k-1&-10k-1\\
0&2k-11&-1\\
0&20k^2-8k-1&-10k-1\\
0&2k-11&-1\\
0&20k^2-8k-1&-10k-1\\
0&2k-11&-1\\
\vdots&\vdots&\vdots\\
0&20k^2-8k-1&-10k-1\\
-1&2k-11&-1\\
0&20k^2-8k-2&-10k-1\\
0&-10k+6&5\\
0&0&1\\
  \end{matrix}
\right)\]
\end{lem}
\bigskip
\bigskip
Also
$\Psi_{2m}=(2m-2).(2m-3).\Psi_{2m-2}.(\overline{2m-3}).(\overline{2m-2})$
and
$\Psi_4 \; = \; \begin{pmatrix}2&1&-1\\7&4&-7\\0&0&\ 1\end{pmatrix}$
so we similarly obtain:

\begin{lem} \label{Psi} For $m\ge 3$, the matrices $\Psi_{2m}$ are given
by
\[\begin{pmatrix} \sp 2&\sp 0&0&\cdots&0&\sp 1&-1\\
\sp 7&-3&0&\cdots&0&\sp 7&-7\\
-7&\sp 4&\UL&{}&\UR&-7&\sp 7\\
\sp 7&-4&{}&{}&{}&\sp 7&-7\\
-7&\sp 4&{}&\parbox[r]{.75in}{\hspace*{.15in}\LARGE{$I_{2m-5}$}}&{}&-7&\sp 7\\
\sp\vdots&\sp\vdots&{}&{}&{}&\sp\vdots&\sp\vdots\\
\sp 7&-4&{}&{}&{}&\sp 7&-7\\
-7&\sp 4&\LL&{}&\LR&-7&\sp 7\\
\sp 7&-4&0&\cdots&0&\sp 8&-7\\
\sp 0&\sp 0&0&\cdots&0&\sp 0&\sp 1
  \end{pmatrix}\]
\end{lem}
\bigskip
\bigskip
Finally, since
\[\Omega_{2m,k}=
\G_{2m-2}^{-k}.\Phi_{2m}.\G_{2m-2}^k.(2m-1).\Psi_{2m}.(\overline{2m-1})\]
we obtain an expression for the monodromy.

\begin{prop}\label{Omega} For $m\ge 3$, $\Omega_{2m,k}$ is given by
\[ \begin{pmatrix}
  a(k)& b(k)&0&\cdots &0 &c(k) &-a(k)+1\\
14k+4 & -8k+1&\UL&{} &\UR &30k+3&-14k-3\\
a(k) & b(k)+1 &{}&{} &{} &c(k) &-a(k)+1\\
14k+4 &
-8k+1&{}&\parbox[r]{.75in}{\hspace*{.07in}\LARGE{$-I_{2m-5}$}}&{}&30k+
3&-14k-3\\
\vdots&\vdots&{}&{}&{}&\vdots&\vdots\\
a(k) & b(k)+1 &{}&{} &{} &c(k) &-a(k)+1\\
14k+4 & -8k+1&\LL&{} &\LR&30k+3&-14k-3\\
a(k) & b(k)+1&0&\cdots &0 &c(k)+1 &-a(k)+1\\
46-70k & -35+40k&0&\cdots &0 &-150k+108 &70k-46\\
7 & -4&0&\cdots &0 &14 &-16
  \end{pmatrix}
\]
where $a(k)=140k^2-64k-10$,
$b(k)=-80k^2+54k+8$, and
$c(k)=300k^2-156k-25$.
\end{prop}

The Alexander polynomial of the 2-component link $L_{2m,k}$ is a function
$\DD_{L_{2m,k}}(t_1,t_2)$ of 2 variables. The reduced Alexander polynomial
is the
single variable polynomial defined by $\Db_{L_{2m,k}}(t)=\DD_{L_{2m,k}}(t,t)$.
For a fibered link $L$ with homology monodromy $\mu$ whose characteristic
polynomial is $p_{\mu}(t)$ one has
\[ \Db_L(t)\cdot (t-1)= p_{\mu}(t) \]
(see \cite{Hill}). For our construction it will suffice to compute
the reduced Alexander polynomials of the
links
$L_{2m,k}$. Again this will be an inductive calculation relying on the
explicit form
for $\Omega_{2m,k}$ given by Proposition~\ref{Omega}.

\begin{thm}\label{DD} The reduced Alexander polynomial $\Db_{L_{2m,k}}(t)$
for the 2-fold covering links $L_{2m,k}$ are given by
\begin{enumerate}
\item[(a)] For $m=2$,\;\; $\Db_{L_{4,k}}(t)= t^2-(140k^2-174k+56)t+1$
\vspace{.1in}\item[(b)] For $m\ge3$,
\begin{multline*}\Db_{L_{2m,k}}(t)=(t^{2m-2}+1)-(140k^2-222k+92)(t^{2m-3}+t)+\\
+(136k^2-258k+119)\sum\limits_{j=1}^{m-2}t^{2j}-
(140k^2-270k+128)\sum_{j=1}^{m-3}t^{2j+1}
\end{multline*}\end{enumerate}
\end{thm}
\begin{proof} Let $p_{m,k}=\det(\Omega_{2m,k}-tI)$, the characteristic
polynomial
of $\Omega_{2m,k}$. First, for $m\ge 4$, we give recursive formulas which
reduce
the calculation of $p_{m,k}$ to that of $p_{8,k}$. According to
Lemma~\ref{Omega},
$\Omega_{2m,k}-tI$ is a rank $2m-1$ matrix of the form
\begin{equation}\Omega_{2m,k}-tI=
\begin{pmatrix}
  a-t& b&0&0&\cdots &0 &c &d\\
x & y-t&-1&0&\cdots&0 &z&w\\
a & b+1 &-t&-1 &{}&0&c &d\\
x & y&0&-t&{}&0&z&w\\
\vdots&\vdots&\vdots&\vdots&{}&\vdots&\vdots&\vdots\\
a & b+1 &0&0 &{} &0&c &d\\
x& y&0&0&{}&-1&z&w\\
a & b+1&0&0&{} &-t &c+1 &d\\
\a & \g&0&0&\cdots &0 &\ve-t &\vt\\
\b & \d&0&0&\cdots &0 &\z &\k-t
  \end{pmatrix}\end{equation}
where we have left out the dependence on $k$. Expand its determinant by the
third
column to obtain
\begin{equation}
           p_{m,k}= \det(U_{m,k})-t\det(V_{m,k})
\end{equation}
Expanding $\det(V_{m,k})$ twice,
each time by the third column, we obtain
\begin{equation}
    \det(V_{m,k})=t^2 \det(V_{m-1,k})
\end{equation}
Similarly, always expanding by the third column, obtain
\begin{eqnarray}
\det(U_{m,k}) = \det(Q_{m,k})-t\det(R_{m,k})\\
\det(Q_{m,k}) = \det(U_{m-1,k})-t\det(S_{m,k})\\
\det(S_{m,k}) = t^2 \det(S_{m-1,k})\\
\det(R_{m,k}) = t^2 \det(R_{m-1,k})
\end{eqnarray}
Equations (2) -- (7) reduce the problem of calculating $p_{m,k}$, $m\ge 4$,
to the
calculation of the quantities in these equations for $m=4$, and this is
accomplished directly from equation (1). Similarly, $p_{6,k}$ is calculated
from
(1) and $p_{4,k}$ is easy to calculate as well. The theorem follows by dividing
$p_{m,k}$ by $t-1$.
\end{proof}

\begin{cor}\label{linking} The two components of the links $L_{2m,k}$ have
nonzero
algebraic linking number.\end{cor}
\begin{proof} The algebraic linking number of the 2 components of $L_{2m,k}$ is
$\Db_{L_{2m,k}}(1)$ \cite{Hill}. This is easily calculated from
Theorem~\ref{DD}:
\[-\Db_{L_{2m,k}}(1)=(4m+132)k^2-(12m+150)k+(9m+36)\]
and the lemma follows simply from this.\end{proof}

\section{Some background on link surgery and Seiberg-Witten
invariants}\label{bg}

The Seiberg-Witten
invariant  of a smooth closed oriented $4$-manifold
$X$ with $b_2 ^+(X)>1$ is an integer-valued function which is defined on
the set of $spin ^{\, c}$ structures over $X$ (cf. \cite{W}). In case
$H_1(X;\Z)$
has no 2-torsion there
is a
natural identification of the $spin ^{\, c}$ structures of
$X$ with the characteristic elements of $H_2(X;\Z)$ (i.e. those elements
$k$ whose
Poincar\'e duals $\hat{k}$ reduce mod~2 to $w_2(X)$).
In this case we view the
Seiberg-Witten invariant as
\[ \ssw_X: \lbrace k\in H_2(X;\Z)|\hat{k}\equiv w_2(TX)\pmod2)\rbrace
\rightarrow \Z. \]
The sign of $\ssw_X$
depends on an orientation of
$H^0(X;\R)\otimes\det H_+^2(X;\R)\otimes \det H^1(X;\R)$. If $\ssw_X(\b)\neq
0$, then $\b$ is called a {\it{basic class}} of $X$. It is a fundamental
fact that the set of
basic classes is
finite. Furthermore, if $\b$ is a basic class, then so is $-\b$ with
$\ssw_X(-\b)=(-1)^{(\text{e}+\text{sign})(X)/4}\,\ssw_X(\b)$ where
$\text{e}(X)$ is
the Euler number and $\text{sign}(X)$ is the signature of $X$.

Now let
$\{\pm \b_1,\dots,\pm \b_n\}$ be the set of nonzero basic classes for $X$.
Consider variables $t_{\b}=\exp(\b)$ for each $\b\in H^2(X;\Z)$ which
satisfy the
relations $t_{\a+\b}=t_{\a}t_{\b}$. We may then view the Seiberg-Witten
invariant
of $X$ as the Laurent polynomial
\[\sw_X = \ssw_X(0)+\sum_{j=1}^n \ssw_X(\b_j)\cdot
(t_{\b_j}+(-1)^{(\text{e}+\text{sign})(X)/4}\, t_{\b_j}^{-1}).\]

  We next recall the link surgery construction of \cite{KL4M}. This
construction
starts with an oriented $n$-component link
$L=\{K_1,\dots,K_n\}$ in $S^3$ and $n$ pairs $(X_i,T_i)$ of smoothly embedded
self-intersection $0$ tori in simply connected $4$-manifolds. The tori are
assumed to be c-embedded, that is, each torus $T_i$ is homologically
essential and
has a pair of embedded curves which generate its first homology and which bound
vanishing cycles (disks of self-intersection $-1$) in $X_i$. For example,
$T$ is
c-embedded if it has a neighborhood $N\C X$ such that the pair $(N,T)$ is
diffeomorphic to $(N_C,F)$ where $N_C$ is a neighborhood of a cusp fiber in an
elliptic surface and $F$ is a smooth elliptic fiber in $N_F$. Let $N(K_i)$ be
disjoint tubular neighborhoods of the components $K_i$ of $L$ in $S^3$ and
$N(L)=\cup N(K_i)$.

Let \,
$\alpha_L:\pi_1(S^3\setminus L)\to \Z $ \,  denote the homomorphism
characterized by the property that it sends the meridian $m_i$ of each
component
$K_i$ to $1$, and
let $\ell_i$ denote the longitude of $K_i$. The
curves $\g_i=\ell_i + \alpha_L(\ell_i) m_i$ on $\bd N(K_i)$
form the boundary of a Seifert surface for the link, and in case $L$
is a
fibered 2-component link, the $\gamma_i$ are given by the boundary
components
of a fiber.

In $S^1\x (S^3\- N(L))$ let $T_{m_i}=S^1\x m_i$,
and define the 4-manifold $X(X_1,\dots X_n;L)$ by
\[
X(X_1,\dots X_n;L)=
(S^1\x (S^3\- N(L))\cup\bigcup\limits_{i=1}^n (X_i\- (T_i\x D^2))
\]
where $S^1\x\bd N(K_i)$ is
identified with $\bd N(T_i)$ so that for each $i$:
\[ [T_{m_i}]=[T_i], \ \ {\text{and}} \ \
[\g_i] = [{\text{pt}}\x\bd D^2]. \]
We have the following calculation of its Seiberg-Witten invariant:

\begin{thm}[\cite{KL4M}] \label{kl} If each $T_i$ is c-embedded in $X_i$ and if
each
$\pi_1(X\setminus T_i) =1$, then
$X(X_1,\dots X_n;L)$ is
simply connected and its Seiberg-Witten invariant is
\[
\sw_{X(X_1,\dots X_n;L)} =
\Ds_L(t_1,\dots,t_n)\cdot\prod_{j=1}^n\sw_{X_j}\cdot(t_j^{1/2}-t_j^{-1/2})
\]
where $t_j=\exp(2[T_j])$ and $\Ds_L(t_1,\dots,t_n)$ is the symmetric
multivariable Alexander polynomial.
\end{thm}

In case each $(X_i,T_i)\cong (X,T)$, a fixed pair, we write \[ X(X_1,\dots
X_n;L)=X_L\] (We implicitly remember $T$, but it is removed from the
notation.)
As an example, consider the case where each $X_i=E(1)$, the rational elliptic
surface ($E(1)\cong \CP\#9\CPb$) and each $T_i=F$ is a smooth elliptic fiber.
Since $SW_{E(1)}=(t^{1/2}-t^{-1/2})^{-1}$, we have that
\begin{equation}\label{E1}
\sw_{E(1)_L}=\Ds_L(t_1,\dots,t_n).
\end{equation}

\section{Symplectic submanifolds}

Let $T$ be a c-embedded symplectic torus in the simply connected symplectic
4-manifold
$(X,\o)$. Then $T$ has a tubular neighborhood which may be identified with
$N=S^1\x
S^1\x D^2$, with $T=S^1\x S^1\x\{0\}$. The symplectic tubular neighborhood
theorem
implies that the restriction of $\o$ to this neighborhood is equivalent to the
symplectic form $dx\wedge dy +rdr\wedge d\theta$.  Let $B$ be a closed
$2m$-strand
braid contained in an unknotted solid torus in $S^3$ with axis $A$.  Define
$T_B$ to
be the torus $T_B=S^1\x B\C N$. Then $T_B$ represents the homology class
$2m\,[T]$.
Furthermore, $T_B$ is a symplectic submanifold of $X$ because its tangent
space at
each point is spanned by $\bd/\bd x$ and the tangent vector $w$ along the
curve $B$,
and because $w$ always has a nontrivial ($\bd/\bd y$)-component,
\[(dx\wedge dy +rdr\wedge d\theta) \,(\bd/\bd x,w)=
dx\wedge dy\,(\bd/\bd x,w)\ne 0.\]

Given $m\ge 2$, consider our family of braids, $B_{2m,k}$ of
\S~\!\ref{braids}. Our
examples are the symplectic tori $T_{B_{2m,k}}$. Let us fix $m\ge 2$ and denote
$B_{2m,k}$ by
$B_k$ and $T_{B_{2m,k}}$ by $\Sig_k$.

It is plausible that one might be able to
distinguish the isotopy classes of the tori $\Sig_k$ by means of the
fundamental
groups of their complements. However, in our situation, where $T$ is
c-embedded, the following lemma points out that one has to work harder.

\begin{lem}\label{pi1} If $T$ is c-embedded then the complement of $\Sig_k$
satisfies
$H_1(X\- \Sig_k;\Z)=\Z_{2m}$ and $\pi_1(X\- \Sig_k)$ is independent of
$k$. If also $\pi_1(X\- T)=1$ then $\pi_1(X\- \Sig_k)=\Z_{2m}$. \end{lem}
\begin{proof} The fundamental group of the complement of $T$ is normally
generated by the boundary $\mu_T$ of a normal disk to $T$. The fundamental
group
of $X\- \Sig_k$ is an amalgamated free product,
\begin{equation}\label{VK}
\pi_1(X\- \Sig_k)= \pi_1(X\- T)*_{\pi_1(\bd(T\x D^2)}\pi_1((T\x D^2)\- \Sig_k)
\end{equation}
and $(T\x D^2)\- \Sig_k$ is the product of a circle with the fiber bundle
$S^1\x _{\vp}\DD$ where $\DD$ is a 2-disk with $2m$ punctures. (Of course,
$\vp$ depends on $k$.) Thus
\begin{equation}\label{FB} \pi_1((T\x D^2)\- \Sig_k)= \la
\mu_1,\dots,\mu_{2m},s,t\; | \; [s,\mu_i]=1,[s,t]=1,
        t\mu_it^{-1}=\vp(\mu_i) \ra
\end{equation}
Since each braid $B_k$ is connected, the action of the monodromy $\vp$ is
transitive on the $\mu_i$. Also, $s$ and $t$ both lie in the image of
$\pi_1(\bd(T\x D^2))$.

If $T$ is c-embedded, then it has a cusp neighborhood $N$, containing both
vanishing cycles, so that the inclusion induces the trivial map
$\pi_1(T)\to\pi_1(\bd N)$. Thus
$s$ and $t$ are both trivial in the amalgamated free product~ \eqref{VK}. It
follows from \eqref{FB} that all the $\mu_i$ are equal in $\pi_1(X\-
\Sig_k)$, and
$\mu_T=\mu_1^{2m}$. The lemma now follows directly from \eqref{VK}.
  \end{proof}
(Note that the most obvious examples, such
as $(X,T) = (E(n),F)$, where $E(n)$ is the elliptic surface over
${\mathbf{CP}}^{\,1}$ without multiple fibers and with holomorphic Euler number
$n$ and $F$ is a smooth fiber, have
$\pi_1(X\- T)=1$.) We shall show that the symplectic tori $\Sig_k$ are not
smoothly isotopic in $X$ by considering
the double branched covers
$\pi_k:\X_k\to X$ branched over the $\Sig_k$.

Let $L_k=\pi^{-1}(A)$ be the double branched covering link. This is the
link that
was denoted $L_{2m,k}$ in \S~\!3.  We may write
\begin{equation}\label{union} \X_k = (X\- N)\cup S^1\x(S^3\- L_k)\cup
(X\- N) \end{equation}
since the double cover is trivial over $X\- N$ and since
$(D^2\x S^1,B_k)= (S^3\- A, B_k)$. In the branched cover \eqref{union}, the
pieces are glued together so that the boundary circles of the fiber
$\Sig_{m-1}\-
(D^2\cup D^2)$ of $S^3\- L_k$ are glued to the boundaries $\bd D^2$ of $X\-
N=X\- (D^2\x T)$. Thus $\tilde{X}_k$ is the manifold $X_{L_k}$ of
\S~\!\ref{bg}. It follows that:
\begin{equation}\label{SWcover} \sw_{\X_k}=\Ds_{L_k}(t_1,t_2)\cdot\sw_{X_1}
\cdot(t_1^{1/2}-t_1^{-1/2})
  \cdot\sw_{X_2}\cdot(t_2^{1/2}-t_2^{-1/2})
\end{equation}
where $T_j$ is a copy of $T$ in the $j$th ($j=1,2$) copy
$X_j$ of $X$, and $t_j=\exp(2[T_j])$.

Assume that there is an isotopy in $X$ which takes $\Sig_i$ to $\Sig_j$. This
isotopy gives rise to a diffeomorphism $h:X\to X$ satisfying $h(\Sig_i)=\Sig_j$
and $h_* = {\text{id}}$ on homology. There is a lift to a
diffeomorphism $\h:\X_i\to\X_j$ of double branched covering spaces.
For a fixed homology class $\b\in H_2(X,\Z)$ consider all the basic
classes $\tilde{\b}$ of
$\X_i$ satisfying
${\pi_i}_*(\tilde{\b})=\b$. For any such
class, it
is also true that ${\pi_j}_*\h_*(\tilde{\b})=\b$ since $h_* = {\text{id}}$. The
invariance of the Seiberg-Witten invariant under diffeomorphisms also implies
that $\ssw_{\X_j}(\h_*(\tilde{\b}))= \ssw_{\X_i}(\tilde{\b})$. Thus for a fixed
$\b\in H_2(X;\Z)$,
\[
\sum_{{\pi_i}_*(\tilde{\b})=\b}\ssw_{\X_i}(\tilde{\b}) =
  \sum_{{\pi_j}_*(\tilde{\g})=\b}\ssw_{\X_j}(\tilde{\g}).
\]

This equation implies that
the Seiberg-Witten invariants of $\X_i$ and $\X_j$ become equal after applying
the projections ${\pi_i}_*$ and ${\pi_j}_*$. Equivalently, working with the
Laurent polynomials as in equation~\eqref{SWcover}, we get
\begin{equation}\label{final}
\Dbs_{L_i}(t)\cdot(\sw_X)^2\cdot (t^{1/2}-t^{-1/2})^2 =
  \Dbs_{L_j}(t)\cdot(\sw_X)^2\cdot (t^{1/2}-t^{-1/2})^2
\end{equation}
where $t=\exp(2[T])$.

\begin{thm} Let $X$ be a symplectic $4$-manifold which contains a c-embedded
symplectic torus $T$. Then each homology class $2m\,[T]$, $m\ge2$, contains the
infinite family $\{ T_{B_{2m,k}}\}$ of symplectic tori, no two of which are
smoothly isotopic.
\end{thm}
\begin{proof} Fix $m$ and consider the double branched covers $\X_k$ of
$(X,T_{B_{2m,k}})$. If $T_{B_{2m,i}}$ is smoothly isotopic to $T_{B_{2m,j}}$,
then equation~\eqref{final} follows.
However, since $X$ is symplectic, $\sw_X\ne0$ \cite{T}, and
it follows from Theorem~\ref{DD} that the $\Dbs_{L_k}(t)$ are all distinct for
different $k$. Thus equation~\eqref{final} can hold only if $i=j$.
\end{proof}

Notice that since we are unable to compute the 2-variable Alexander
polynomials for the
links $L_i$, this proof, in itself, does not show that the covers,
$\X_k$ are mutually nondiffeomorphic --- only that they can not be
made diffeomorphic
via a
$\Z_2$-equivariant diffeomorphism which covers the identity on $H_2(X;\Z)$.
Thus, without further information about the Alexander polynomials of the links,
we are unable to show via this technique that there is no diffeomorphism of $X$
which throws $T_{B_{2m,i}}$ onto $T_{B_{2m,j}}$ for $i\ne j$. However for the
case $X=E(1)$, we get a stronger result.

\begin{thm} \label{E1thm} Let $T$ denote a smooth elliptic fiber in the
rational elliptic surface, $E(1)$. Then each homology class $2m\,[T]$, $m\ge2$,
contains the infinite family $\{ T_{B_{2m,k}}\}$ of symplectic tori, no two of
which are equivalent under diffeomorphisms of $E(1)$.
\end{thm}
\begin{proof} In this case the double branched covers $\X_k$ cannot be
diffeomorphic for different $k$. For, it follows from equations~\eqref{SWcover}
and \eqref{E1}
that $\sw_{\X_k}=\Ds_{L_k}(t_1,t_2)$. Thus $\Ds_{L_k}(1,1)$ is a diffeomorphism
invariant of $\X_k$, the sum of all its Seiberg-Witten invariants. The
calculation
of Corollary~\ref{linking} shows that these numbers are different for different
$k$.
\end{proof}

\section{Lefschetz fibrations}

In this section, we show how
our constructions above naturally yield examples of infinite classes of
homeomorphic but
nondiffeomorphic symplectic manifolds, all of which admit Lefschetz
fibrations of
fixed fiber genus. (There is a more general construction presented in
\cite{FSLef}.) For simplicity, we
restrict ourselves with the
application of this
procedure to the rational elliptic surface $E(1)$. Let $T=F$, a generic
elliptic
fiber in $E(1)$ and let $\tilde{X}_{2m,k}$ be the double branched cover of
$E(1)$
with  branch set $T_{B_{2m,k}}$. Then $\tilde{X}_{2m,k}\cong
E(1)_{L_{2m,k}}$ is a homotopy $K3$ surface.

It is well-known that $E(1)$ admits a genus $0$ fibration with 4 singular
fibers.
This is seen by noting that $E(1)$ is the double branched cover of $\SS$
with branch
set equal to
$4$ disjoint copies of $S^2\x \{{\rm{pt}}\}$ together with $2$ disjoint
copies of
$\{{\rm{pt}}\}\x
S^2$. The resultant branched cover has $8$ singular points (corresponding
to the
double points in the branch set), whose neighborhoods are cones on
${{\mathbf{RP}}^{\,3}}$. These are desingularized in the usual way,
replacing these
neighborhoods with cotangent bundles of $S^2$. The result is $E(1)$. The
horizontal
and vertical fibrations of
$\SS$ pull back to give fibrations of
$E(1)$ over
${{\mathbf{CP}}^{\,1}}$. A generic fiber of the vertical fibration is the
double
cover of
$S^2$, branched over $4$ points --- this gives an elliptic fibration on
$E(1)$. The
generic fiber of the horizontal fibration is the double cover of
$S^2$, branched over $2$ points --- this gives the genus $0$ fibration of
$E(1)$.
The $4$ singular fibers are the preimages of the four $S^2\x
\{{\rm{pt}}\}$'s in
the branch
set. The generic fiber $T$ of the elliptic fibration meets a generic fiber
$\Sig_0$
of the horizontal fibration in $2$ points, $\Sigma_0\cdot T=2$. This means that
$T_{B_{2m,k}}$ meets $\Sig_0$ transversely in $\Sig_0\cdot T_{B_{2m,k}}
=4m$ points.
Therefore the horizontal fibration on $\SS$ lifts to a fibration
on $\tilde{X}_{2m,k}$ whose generic fiber is the double cover of $S^2$
branched over
$4m$ points, that is, a genus $2m-1$ fibration. The definition of a Lefschetz
fibration requires the monodromy around each singular fiber to be a Dehn
twist. This
is not true for these examples, but they can be perturbed to be Lefschetz (see
\cite{GS}).

One can give an alternative description of this fibration on
$\tilde{X}_{2m,k}\cong
E(1)_{L_{2m,k}}$. The elliptic fiber $T$ of $E(1)$ meets each genus~ $0$ fiber
transversely, and it meets a generic genus~ $0$ fiber $\Sig_0$ twice. Thus
$\Sig_0\- (\Sig_0\cap T)$ is an annulus. The construction of
\[
E(1)_{L_{2m,k}}
=\{ E(1)\- N(T)\}\,\cup \{S^1\x M_{L_{2m,k}}\- (N(T_1)\cup
N(T_2))\}\,\cup\{E(1)\-
N(T)\}
\]
preserves the fibrations. In the manifold (with boundary)
\[\{ E(1)\- N(T)\}\,\cup \{S^1\x M_{L_{2m,k}}\- (N(T_1)\cup N(T_2))\}\]
a generic fiber of the induced fibration is the union
of 2 fibers of the fibration
\[\begin{CD}
\Sig_{m-1}'' \to\; &S^1\x M_{L_{2m,k}}\- (N(T_1)\cup N(T_2))&\; =
            S^1\x (S^3\- L_{2m,k})\\
  & @VVV&\\
& S^1\x S^1&\\
\end{CD}
\]
together with the annulus $\Sig_0\- (\Sig_0\cap T)$. This is a surface of genus
$2m-2$ with 2 boundary components. Adding the second copy of $E(1)$ adds an
annulus which closes up the surface,
and we obtain a surface of
genus $2m-1$.

It is not difficult to see that the fibrations on $\tilde{X}_{2m,k}$ which
are described
here are actually hyperelliptic, that is, the hyperelliptic involutions on the
fibers extend to a global involution of $\tilde{X}_{2m,k}$. In fact, the
orbit space of this hyperelliptic
involution is $E(1)$ and the image of the fixed point set is just the torus
$T_{B_{2m,k}}$.

\end{document}